# Fast Sparse Matrix-Vector Multiplication on GPUs: Implications for Graph Mining

Xintian Yang, Srinivasan Parthasarathy, P. Sadayappan
Department of Computer Science and Engineering
Ohio State University, Columbus, OH 43210
{yangxin, srini, saday}@cse.ohio-state.edu

## ABSTRACT

Scaling up the sparse matrix-vector multiplication kernel on modern Graphics Processing Units (GPU) has been at the heart of numerous studies in both academia and industry. In this article we present a novel non-parametric, self-tunable, approach to data representation for computing this kernel, particularly targeting sparse matrices representing power-law graphs. Using real web graph data, we show how our representation scheme, coupled with a novel tiling algorithm, can yield significant benefits over the current state of the art GPU efforts on a number of core data mining algorithms such as PageRank, HITS and Random Walk with Restart.

## 1. INTRODUCTION

Over the last decade we have witnessed a revolutionary change in the way commodity processor architectures are being designed and implemented. CPUs with superscalar out-of-order execution, vector processing capabilities, and simultaneous multithreading, chip multiprocessing (CMP), and high-end graphics processor units (GPU) have all entered the mainstream commodity market. Data intensive algorithms often require significant computational resources, and thus stand to benefit significantly from such innovations if appropriately leveraged.

In this article we develop a novel approach to facilitate the efficient processing of key graph-based data mining algorithms such as PageRank [15], HITS [10] and Random Walk with Restart[18] on modern GPUs. A common feature of these algorithms is that they rely on a core sparse matrix vector multiplication kernel (SpMV). Implementations of this kernel on GPUs have received much attention recently from the broader scientific and high performance computing communities [5, 2, 6] including an industrial strength effort from NVIDIA research [3].

The key difference between past work and ours is that here we are interested in the processing of sparse matrices that represent large graphs – typically with power-law [12] characteristics. This difference is also central to the specific architecture-conscious approach we propose for processing the SpMV kernel. We transform and represent the matrix in such a way so as to facilitate tiling – a key strategy used to enhance temporal locality. Additionally we rely on a composite storage algorithm that leverages the skew in the degree distribution afforded by the fact that these matrices represent power-law graphs. Architectural features of the GPU such as the texture cache are also effectively leveraged in the processing of the kernel. Moreover, we also demonstrate how the basic approach can be extended to handle web-scale graphs with billions of edges that do not fit in the memory of a single GPU by suitably leveraging multiple GPUs.

The method summarized above relies on significant programmer expertise to tune two key parameters of the approach. To alleviate this we present a systematic mechanism for tuning the parameters automatically at runtime depending on input matrix characteristics. Unlike existing work on parameter auto-tuning for SpMV kernel [6, 19], an important side effect of our approach is a reliable performance model for predicting overall performance of the kernel.

We present a comprehensive empirical evaluation of the proposed methods on three data mining algorithms and the base SpMV kernel on a range of real datasets including several web-scale graph datasets. We find that on moderately sized datasets that fit on a single GPU, our SpMV kernel as well as the graph mining methods that rely on this kernel – HITS, PageRank and Random Walk with Restart – are typically **1.8 to 2.1** times faster than an industrial strength state-of-the-art GPU competitor and anywhere from **18 to 32** times faster than a similarly structured and optimized CPU implementation. We find that our methodology can *scale quite comfortably to webscale datasets with billions of edges*, on multiple GPUs with parallel efficiencies of up to **70%**. Finally, we empirically demonstrate the effectiveness of our autotuning approach both in terms of its *ability to correctly select parameters for our kernel* on a wide range of datasets, and in terms of its ability to *reliably predict the absolute performance of the kernel* under different parametric settings, suggesting it can be used for adaptive algorithm designs in next generation hybrid architectures[16].

## 2. RELATED WORK

Several existing efforts have targeted optimizing the SpMV kernel for the GPU [3, 2, 6]. However, none of the above takes into account the skew of the non-zero distribution present in matrices representing power-law graphs. The per-formance of previously reported implementations is low on such matrices due to their power-law characteristics.

Bell and Garland [3] propose representations of sparse matrices on the CUDA platform in NVIDIA's SpMV library. Their library includes CSR, CSR-vector, COO, ELL, HYB, DIA and PKT formats, and is the most closely related work to our research (see Appendix B). Choi [6] proposes a *blocked ell-pack* format in which the non-zeros are stored in blocks and the blocks are indexed with ELL format. This introduces high memory overhead on power-law matrices.

Because of the importance of the SpMV kernel, researchers have put substantial efforts on optimization techniques over various platforms. Vuduc *et al* [19] study optimizations and performance auto-tuning in single core CPUs. Williams [20] compares SpMV kernels on emerging multicore platforms.

Previous studies emphasize performance optimization of the SpMV kernel on various platforms. However, they require parameter tuning to achieve high performance, which depends on the characteristics of the input matrices. Therefore, an important issue is the development of a performance model to guide the tuning. Hong [8] proposed an analytical model of GPU architecture. An adaptive performance modeling tool based on the work flow graph of a GPU application is proposed by Baghsorkhi [1]. Choi [6] proposed a model-driven autotuning framework in their SpMV work.

## 3. METHODOLOGY

In this section we first present optimization techniques for SpMV kernel on power-law matrices within a single GPU. Next, we show how our SpMV kernel can be extended to handle out-of-core matrices on a multi-GPU cluster. Finally, we present a systematic mechanism for automatically tuning the parameters in our methods. We assume that readers are familiar with the CUDA architecture (Appendix A).

### 3.1 Single-GPU SpMV

Our optimizations are based on a series of observations from benchmarks that demonstrate the limitations of previous work. We propose solutions that target these limitations and thereby improve the performance.

The SpMV kernel is a bandwidth limited problem since the floating point operations per memory access is low. When computing the product of a sparse matrix $A$ and a vector $x$, the memory accesses to matrix $A$ are optimized to be fully coalesced in NVIDIA's SpMV library. But the accesses to vector $x$ have never been optimized. Also vector $x$ is the only reusable memory in the SpMV kernel.

**Observation 1: Each row accesses random elements in vector $x$.** In the adjacency matrix of a power-law graph, the column indices of the non-zeros on each row are not continuous, and are relatively random, which leads to non-coalesced memory addresses when accessing $x$. Previous work [3, 2] has bound the entire vector $x$ to the texture memory and utilizes the cache of texture memory to improve the locality. But the size of $x$ is usually much larger than the size of the texture cache. The resultant cache misses reduce memory bandwidth utilization due to the long latency of non-coalesced global memory accesses.

**Solution 1: Tiling matrix $A$ and vector $x$ with texture cache.** Suppose we divide matrix $A$ into fixed width tiles by column index and segment vector $x$ correspondingly, so that each tile of $A$ only needs to access one segment of $x$.

If one segment of $x$ can fit in the texture cache, once the elements are fetched into the cache, none of them will be kicked out until the computation of this tile finishes. Therefore, we can get maximum reuse of $x$.

The texture cache size is a key factor in determining the width of the tiles. To estimate the texture cache size (since this is not provided by the manufacturer) on our Tesla GPU, we conduct benchmarking experiments as follows. We *mod* the column indices of a large sparse matrix by tile width, so all accesses to vector $x$ are mapped to one tile. We vary the tile width from $100K$ to $1K$ and run the multiplication. The performance improves most significantly when tile width = $64K$, corresponding to 256 KB of cache size. So our tile width is fixed to $64K$ columns.

The performance of tiling the entire matrix $A$ and vector $x$ is still low. Since we divide all the columns of matrix $A$ into tiles, there could be too many tiles when the matrix is large. Each tile needs to add its partial result to the final result $y$, leading to non-coalesced memory accesses overhead. Also, the write-back result of one tile has to be visible to the next tile before the next tile can start; otherwise, memory read-after-write conflicts could happen. To avoid such memory conflicts, we restart a kernel for each tile, which also causes an overhead. Therefore, tiling the columns of the matrix $A$ fails to improve the performance of the SpMV kernel.

**Observation 2: Column lengths follow a power-law distribution.** Suppose a matrix is the adjacency matrix of a power-law graph, the number of non-zeros in the columns of the matrix will follow a power-law distribution. So there are large number of columns with few non-zeros. In tiles containing such columns, we cannot get much reuse of vector $x$, but we still need to restart large number of kernels to compute them and incur the overheads. On the other hand, the long columns in the power-law distribution concentrates large portion of the non-zeros in the matrix. There is a lot of memory reuse of vector $x$ in such columns. The denser the columns, the more the benefits. If we can first tile such columns greedily, we can finish the majority of the total computation efficiently in the first few tiles. The overall performance of the entire matrix will be improved.

**Solution 2: Reorder columns by column lengths and partially tile $A$.** Our idea is to first reorder the matrix columns in decreasing order of length. We can divide the reordered matrix into two sub-matrices by setting a threshold on column length. Long columns form a denser sub-matrix; the remaining short columns form a sparser sub-matrix. The denser sub-matrix contains a lot more non-zero elements and fewer columns than the sparser sub-matrix. According to Amdahl's law, the overall performance of the SpMV kernel will be improved if the computation in the denser sub-matrix can be finished efficiently. Now we can tile the denser sub-matrix with texture cache. The non-zero elements are concentrated in a small number of tiles so that we can still gain the benefits from $x$ vector caching as well as avoid the overhead of initializing too many tiles. Here we introduce a threshold parameter to partition the matrix into two sub-matrices. We will discuss how to automatically determine this parameter based on the distribution of the column lengths of the matrix in Section 3.3 .

Figures 1(a) to 1(c) illustrate the above transformation procedure on a small sparse matrix. Figure 1(a) is the original matrix; Figure 1(b) reorders the columns of the matrix in decreasing order of column lengths. In this example, we set

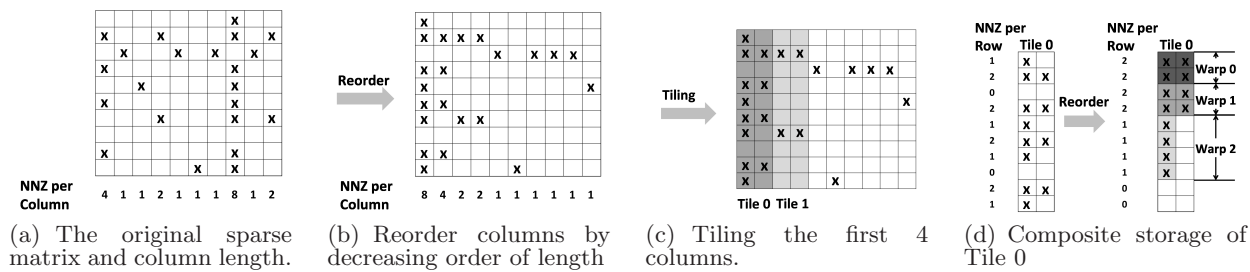

(a) The original sparse matrix and column length.
(b) Reorder columns by decreasing order of length
(c) Tiling the first 4 columns.
(d) Composite storage of Tile 0

Figure 1: Illustrative example of tiling and composite storage. X − non-zeros

the column length threshold to 2. Columns with more than or equal to 2 non-zero elements will be placed in the denser sub-matrix; the other columns with only 1 non-zero element will be placed together in the sparser sub-matrix. Suppose the texture cache can only hold 2 floating point numbers in this example, the denser sub-matrix with 4 columns will be partitioned into 2 tiles as shown in Figure 1(c).

Amongst all the kernels in NVIDIA's SpMV library, HYB and COO perform best on power-law matrices. The computation in the sparser matrix is run under the HYB kernel, because HYB has the best performance. The computation within each tile of the denser matrix will be performed using COO kernel. The resulting vector $y$ from the denser and sparser sub-matrices will be combined to the final result.

**Observation 3: Performance of COO kernel is limited by thread divergence and serialization.** When computing each tile, the COO kernel cannot utilize the massive thread level parallelism in CUDA efficiently although it is more efficient than the CSR-vector and ELL kernel on such data. In the COO kernel, the inputs are three arrays storing the row indices, the column indices and the values of non-zero elements in the matrix. These three arrays are all divided into equal length intervals. Each interval is assigned to one warp. Note that this partition only equally distributes workload to warps. It does not consider that a row may cross the boundary between two warps. Each warp of threads iterates over an interval in a strided fashion. The stride size equals warp size, and a thread within a warp only works on one element in one stride. A thread first fetches the value in the $x$ vector based on the column index, and then multiplies the $x$ value with the element in matrix $A$ and stores the result in a shared memory space reserved for it. The next step is the *sum reduction* of the multiplication results within one stride by a binary reduction operation. But one stride can contain non-zeros from more than one row. When the reduction operation tries to add two operands, it has to first check whether the two operands are from the same row in the original matrix. If not, this warp of threads will be serialized due to the thread divergence. This leads to low thread level parallelism in the COO kernel.

**Observation 4: Performances of CSR-vector and ELL kernel are limited by imbalanced workload.** The CSR-vector kernel performs the best when the rows of a matrix are long and with similar length. Non-zeros are stored in row major in CSR format. CSR-vector kernel assigns one warp of threads per row. Each warp iterates on the row with stride size the same as warp size, and performs multiplication and summation operations. After the last iteration on this row, the threads in a warp perform a binary reduction to obtain the final result of this row. In all the summation and reduction operations, the threads do not need to check whether the two operands are from the same row. However, CSR-vector kernel is most efficient when the number of non-zeros in a row is an integer multiple of the warp-size.

The ELL kernel achieves peak performance if there are large number of short rows with similar lengths in the sparse matrix. In the ELL format, all rows have the same length and 0s are padded to the rows shorter than this length. The non-zeros are stored in column major order. A warp of threads is assigned to work on 32 consecutive rows, with each thread working on the multiplication and reduction of one row. The threads within a warp iterate over the columns efficiently with hardware synchronization.

**Observation 5: Tile row lengths follow power-law.** Due to the scale-free property of power-law graphs, we observe that after tiling, the row length within a tile also follows a power-law distribution. We propose to address this fact via a novel storage format of matrix $A$ within each tile to further improve the efficiency of SpMV kernel.

**Solution 3: Composite tile storage scheme.** Our composite storage scheme combines the CSR and ELL formats as follows. Our algorithm starts by ranking the row lengths from high to low. A workload size is defined as the total number of non-zeros in the longest row or several long rows at the top of the ranking, depending on the dataset. We will discuss the auto-tuning of this parameter in Section 3.3. Then rows in a tile will be partitioned into approximately balanced workloads. This can be implemented by traversing the row length ranking from top to bottom. A new row is packed into a workload until it exceeds the workload size, then a new workload is initialized. Each workload can be viewed as a rectangle, where the width $w$ is the length of the first row (the longest row in this workload) and the height $h$ is the number of rows in this workload. If $w \geq h$, this workload will be stored in row major in global memory and computed by CSR-vector kernel; otherwise, it will be stored in column major and computed by ELL kernel. Note that if a workload is stored in row major, all rows will be padded to the same length as $w$ with 0s; and 0s will also be padded to ensure that $w$ (or $h$) is an integer multiple of warp size when a workload is stored in row (or column) major. After the above partition and transformation of storage format, each workload is assigned to a warp of threads and computed with the most suitable kernel.

The above composite storage scheme is designed for tiles of the denser sub-matrix in the original matrix. However, it can be used for all matrices whose row lengths follow a power-law distribution. We observed that the row lengths in the sparser sub-matrix also follow a power-law distribution. Therefore, we also transform the sparser sub-matrix as one matrix tile into the composite storage format.

Figure 1(d) illustrates how tile 0 from Figure 1(c) is transformed to our composite storage on a fictitious architecture with two threads per warp. The rows in tile 0 are first reordered by row length. Suppose we set the workload size to be 4. The first two rows are packed into the first workload, stored in row major and assigned to warp 0 for computation. The two threads in warp 0 first do multiplication and reduction on row 0 using CSR-vector kernel and move to row 1 together. The next two rows are packed together, stored and computed by warp 1 in a similar fashion. The remaining four single element rows are stored in column major and computed by warp 2. The two threads in warp 2 start from the first two rows vertically using ELL kernel and then move to the last two rows.

**Elimination of Partition Camping:** The global memory is divided into 8 memory partitions of 256-byte width. Concurrent memory requests to the global memory should be distributed uniformly amongst partitions. The term *partition camping* [14] is used to describe the situation when global memory accesses are congested and queued up at some partitions while the other partitions are idle. All data in strides of 2048 bytes (or 512 floats) map to the same partition. In our tile-composite format, if the workload size is an integer multiple of 512 floats, then the start memory addresses of all workloads are mapped to the same partition. All warps will queue up at each partition when the warps iterate on their own workload. Thus, we will have the problem of partition camping. To avoid this problem, we add 256 bytes of memory to the end of each workload if the workload size is an integer multiple of 512 floats.

**Sorting Cost:** Sorting is used to re-structure the matrix to improve memory locality. The cost of sorting is relatively cheap when the rows and columns follow power-law. The rows and columns can be bounded by some small number $k$ in the long tail of the power law distribution. These rows or columns can be sorted by counting sort in linear time [7]. The remaining rows or columns can be sorted very quickly. Moreover, we only need to perform the sorting once as a data preprocessing step. In applications such as the power method where the SpMV kernel is called iteratively until the result converges, the cost of sorting can be amortized.

### 3.2 Multi-GPU SpMV

In this section, we show how our SpMV kernels can scale to large web graphs that cannot fit in the memory of a single GPU. To handle out-of-core matrices, we can either use a single GPU to work on chunks of the matrix in serial, or distribute the chunks to multiple GPUs. Because the single GPU strategy has to move the data from CPU to GPU in every iteration, the bandwidth of the PCI-Express bus from CPU to GPU (8GB/s) will become the performance bottleneck in the single GPU kernel, because our best kernel can comfortably achieve 40GB/s bandwidth (see Figure 2(b)). Hence we devise a method that can use a multi-GPU cluster to compute SpMV kernel on large-scale dataset.

In the cluster, each node keeps one local partition of the matrix. At the end of each iteration, all nodes need to broadcast their local result vector $y$ to the other nodes, so that the other nodes can update their local copy of the vector $x$ for the next iteration. The communication cost of broadcasting, decided by the partition scheme, is the key factor that limits the scalability of the multi-GPU SpMV kernel. Algorithms such as matrix partitioning and graph clustering can be used to minimize this cost. But those algorithms are often more expensive than the iterative SpMV kernel which would be self-defeating. Hence we only consider simple partition schemes such as partition by rows, by columns and by grids.

The communication cost is lower if the *matrix is partitioned by rows* rather than by columns. Suppose we have N rows and P processors. If the matrix is partitioned by rows, each processor only needs to send out N/P elements of vector $x$. But if partitioned by columns, all processors need to send out N elements. Also partitioning by rows does not necessitate any reduction operations after vector $x$ is gathered. For similar reasons, we can show that partitioning by rows is superior to partitioning by grids. In our multi-GPU SpMV kernel, we choose to partition the matrix by rows with a partition scheme that can assign approximately equal number of rows and equal number of non-zeros in each partition. Such scheme can guarantee both balanced computation workload (number of non-zeros) and balanced communication cost (number of rows) on each node. We use bitonic partitioning [17] for this task. The intuition of bitonic partitioning is as follows: The matrix rows are first sorted by length. Each iteration of the algorithm processes P rows and assigns them to P processors. The processor that got the longest row in the previous iteration will get the shortest row in the current iteration.

Any SpMV kernel can be plugged into this multi-GPU framework to perform local computation. Because power-law graphs are known to be scale-free, we observe that the rows and columns of each partition of a power-law matrix also follow power-law. Hence we can expect our optimized SpMV kernel for power-law matrices to be a good fit for the local computations at each GPU in the cluster.

We next introduce an automatic parameter tuning method to find the optimal parameters for each local kernel, since it is prohibitive to use exhaustive search to find optimal parameter settings for the kernels on every node.

### 3.3 Automatic Parameter Tuning

The practical utility of our optimization approach is somewhat limited by the need to carefully tune two parameters – the number of tiles and the partitioning strategy of each tile. We address this limitation as follows.

The first parameter we need to determine is the number of tiles in the composite kernel. It turns out that this parameter is relatively straightforward to estimate. The performance gain of the tiling strategy hinges on the reuse that comes from the temporal locality of the texture cache when accessing vector $x$. If a matrix column has only one non-zero element, there is no reuse benefit. It turns out that if there is any reuse, however small, tiling is typically beneficial. Thus, a new tile should not be added if its first column has only a single element (line 7–8 Algorithm 1 in Appendix E).

Next we need to determine how to partition each tile into small rectangular workloads and assign each workload to a warp of threads in the SpMV kernel. Since the tile-composite kernel requires a balanced workload for each warp, the size of one workload directly decides the partition of the tile. The *problem of partitioning can be transformed to finding the optimal workload size*. Considering the search space of the workload size, the lower bound is the length of the first row in the tile. This is because the rows within each tile are reordered in decreasing order of row lengths, and the longest row cannot be partitioned into two workloads.

The upper bound of the workload size is the total number of non-zeros in the tile divided by the maximum number of active warps available (960 on the Tesla GPU) to fully utilize the GPU resources. Apart from the upper and lower bounds, an additional constraint is that the workload size must be an integral multiple of the first row in this tile (line 11 Algorithm 2 in Appendix E). This is because the first workload of this tile must be a rectangle area of non-zeros where each row is padded to the same length as the first row. These constraints dramatically reduce the search space of viable workload size settings. We must now estimate the performance under each setting and return the one that is predicted to perform the best (Algorithm 2 in Appendix E). We next describe our performance model which has broader utility beyond just our use for parameter setting.

**Performance Model:** *Our model relies on an execution model of CUDA kernels - which forms an offline component in our model - and the non-zero element distribution of the input matrix - which is an online component.* The offline component seeks to create a lookup table indexed on the rectangular shape ($w$ columns, $h$ rows) of a workload and its corresponding performance. At runtime the online component given a particular input matrix tile computes an estimate of its runtime cost for different viable workload sizes. The lowest cost option is selected as the partitioning strategy. We next detail both steps.

Given a rectangle workload whose shape is defined by $w$ and $h$, we construct a lookup table establishing a mapping from the shape of the workload to its performance on one thread warp. We run offline benchmarks to establish this mapping as follows. Each benchmark is for one combination of $w$ and $h$, and we artificially construct a matrix in tile-composite format, in which all workloads are set to the same $w$ by $h$ shape and there are large number of such workloads to fill the computation pipeline and to hide the memory latency. We measure the performance of all realizable combinations of $w$ and $h$. This may seem to be an exhaustive process but fortunately we have several constraints that limit the combinations we need to evaluate. Moreover it is a one time offline cost that is independent of the dataset. The first constraint is the upper bound of the workload size. In practice, we can choose a large enough upper bound to cover all possible workloads. The second constraint is that either $h$ or $w$ must be a multiple of 32. Therefore the total number of combinations of $w$ and $h$ is relatively small and finite. To reiterate, we only need to construct this mapping once for a given GPU architecture (e.g. Tesla) and this mapping can be repeatedly used for any input data matrix.

In the online component of our model, given a particular input matrix tile we first identify viable workload sizes based on the upper and lower bounds as noted above. Then for a specific workload size we partition the tile into workloads of potentially different shapes but of roughly the same size, in the same way as the tile-composite kernel (Figure 1(d), line 8 – 9 Algorithm 3 in Appendix E). We then estimate the performance of the entire tile as the average of the performance of all its constituent workloads as follows.

After we establish the mapping from $w$ and $h$ to a performance number, we simulate the computation of a matrix tile by looking up performance numbers from the mapping instead of running the actual tile-composite kernel. Here we need to consider the execution model of CUDA hardware. Each streaming multiprocessor (SM) can only serve $MAX\_ACT\_WARP/SM$ number of active warps at a time. So if the total number of warps is larger than the maximum number of active warps on a GPU, the thread warps will be divided into iterations so that each iteration can fit in all the SMs on the GPU. The number of iterations can be computed as (line 5 Algorithm 3 in Appendix E):

$$I = \left\lceil \frac{TOTAL\_WARP}{MAX\_ACT\_WARP/SM * NUM\_SM} \right\rceil \quad (1)$$

The total running time will be the sum of the running time in each iteration (line 19 Algorithm 3 in Appendix E):

$$t = \sum_{i=1}^{I} t_i \quad (2)$$

The running time of each iteration $i$ is computed as the total size of the workload in this iteration divided by the average performance of all the warps in this iteration (line 18 Algorithm 3 in Appendix E):

$$t_i = Size(i)/P_i \quad (3)$$

where $Size(i)$ is the total size of all the workloads in iteration $i$ (line 13 Algorithm 3 in Appendix E),

$$Size(i) = \sum_{warp(j) \in Iter(i)} w_j * h_j \quad (4)$$

and $P_i$ is estimated by the average performance of all warps in this iteration (line 18 Algorithm 3 in Appendix E),

$$P_i = \frac{\sum_{warp(j) \in Iter(i)} Performance(w_j, h_j)}{Total\_Warps(i)} \quad (5)$$

For given $w_j$ and $h_j$, $Performance(w_j, h_j)$ can be found in the mapping we established from the offline benchmarks.

Based on Equation $1 - 5$, we can calculate the prediction of the total running time of computing a matrix tile with the given workload size. Our parameter auto-tuning method will call the above performance model as a subroutine to estimate the optimal setting of workload size in each tile (line 7 Algorithm 2 in Appendix E).

Note that the thread block size and the number of thread blocks are parameters in all CUDA kernels. In our tile-composite kernel, we assume full occupancy of all multiprocessors, which means there are 32 active warps on each SM at any time. This assumption is reasonable because if there are fewer warps on each SM, the performance will decrease. The maximum size of a thread block is 512 (=16 warps), and each SM can have 8 active thread blocks at maximum. So the 32 active warps on each SM can only be organized as 8 blocks with 4 warps per block, 4 blocks with 8 warps per block, or 2 blocks with 16 warps per block. In practice, the three configurations yield similar performance because the physical executions of the warps are the same on the hardware. In our experiments, we choose 8 warps per block.

The above performance modeling method is designed for the tiling part of the tile-composite kernel. A similar method is used to model the sparse part of the matrix by running offline benchmarks without using the texture cache.

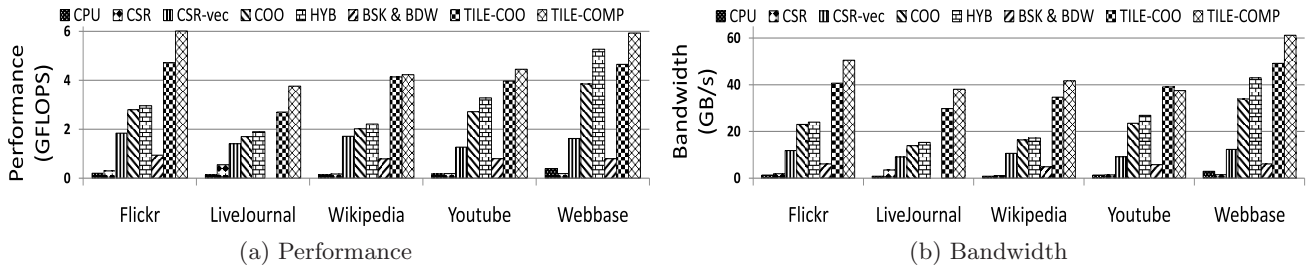

Figure 2: SpMV kernels comparison on matrices representing power-law graphs.

## 4. EXPERIMENTS

In this section, we present experimental results of our optimizations on SpMV kernels and auto-tuning method. The details of the datasets and the hardware configurations are introduced in Appendix C.

### 4.1 Single-GPU SpMV Kernel

We compare a CPU implementation of the CSR kernel, all six GPU kernels from NVIDIA's SpMV library, Baskaran and Bordawekar's optimized CSR GPU kernel (BSK & BDW) and our two optimized GPU kernels (TILE-COO and TILE-COMPOSITE) on the matrix datasets in Table 2. We report the speed of execution in GFLOPS and bandwidth in GB/s. The running time is averaged over 500 iterations. All kernels are run in single precision. Binding the entire vector $x$ to texture cache performs consistently better than not binding in all NVIDIA's kernels [3] and BSK & BDW's kernel [2]. So we only report the performance of these kernels with texture cache binding. We use 256 threads per thread block. This setting is default in NVIDIA's SpMV library. Under this setting, there are enough warps in each thread block to hide the memory latency. In our tiling method, we have to set a threshold to decide the number of tiles. In this section, we choose this parameter by exhaustive search to present the best performance that can be achieved by our kernel. Later, we will examine the performance under auto-tuning. The performance of the SpMV kernels on power-law matrices are shown in Figure 2. The results on the other matrices are presented in Figure 7 Appendix D.

**Performance on power-law matrices**: Our tile-coo and tile-composite methods clearly dominate the other kernels on the Flickr, LiveJournal and Wikipedia datasets. Our tile-composite kernel has an average **1.95x** speedup over NVIDIA's best kernel – HYB kernel on these datasets. On Webbase and Youtube datasets, which are small power-law matrices, COO and HYB kernel perform close to our optimizations. Our Tile-composite kernel is about **13%** faster on Webbase and **36%** faster on Youtube than the NVIDIA HYB kernel. From Table 2, we can see the numbers of rows and columns are low in the Webbase and Youtube matrices, and also the numbers of non-zeros per row and column are low. These properties of the Webbase and Youtube matrices hide the advantages of our optimizations for the following reasons. First, there is little reuse of vector $x$ if non-zeros per column is low. This leads to lesser benefit from our tiling optimization. Second, when the number of columns is small, COO and HYB kernel have better probability of cache hits when they bind the entire vector $x$ to the texture cache. Third, the total number of non-zeros in a tile is low so our composite storage scheme will pad more zeros and thereby cause memory access overhead. We do not report perfor-mance of NVIDIA's PKT kernel on these datasets since the partition step within this kernel does not produce balanced enough *packets* and leads to kernel failure.

| Graph | CPU | COO | HYB | TILE-COO | TILE-Comp |
|---|---|---|---|---|---|
| Flickr | 23.99 | 1.67 | 1.60 | 0.90 | 0.83 |
| LiveJournal | 82.23 | 6.19 | 5.57 | 3.75 | 3.44 |
| Wikipedia | 52.12 | 2.99 | 2.83 | 1.76 | 1.63 |
| Youtube | 11.81 | 0.72 | 0.66 | 0.68 | 0.65 |

Table 1: Running time of PageRank (in seconds)

### 4.2 Graph Mining Applications

In this section we discuss how the single GPU performance of our SpMV kernel translates to performance on important graph mining algorithms. Due to space limitations, in this section we present results for PageRank, similar results for Random Walk with Restart (RWR) and HITS are presented along with implementation details for all three algorithms in Appendix F. We implement these algorithms using four SpMV kernels: COO, HYB, TILE-COO, and TILE-Composite kernels. These four kernels are generally the top performing kernels from the experimental results presented in the previous section.

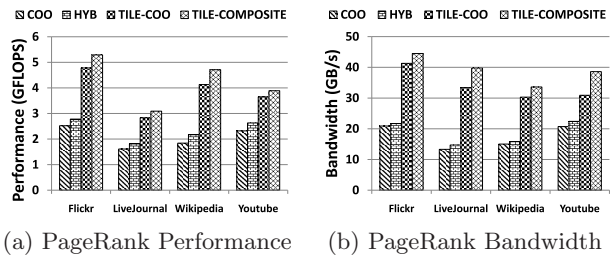

(a) PageRank Performance  (b) PageRank Bandwidth

Figure 3: Performance and bandwidth of PageRank.

We run Equation 6 of Appendix F iteratively with the SpMV kernels and check whether $p$ converges on the 4 graph datasets in Table 2. The speed and bandwidth performance of PageRank based on the four kernels are shown in Figure 3(a) and Figure 3(b). The total running time on each graph is shown in Table 1 in comparison to a CPU implementation of PageRank. Our tile-coo and tile-composite kernel achieve about **2x** speedup over COO and HYB kernel on Flickr, LiveJournal and Wikipedia graphs. Our optimizations are marginally better than NVIDIA's COO and HYB kernel on the Youtube graph (reasons noted earlier). Compared with the CPU PageRank, all GPU implementations achieve between **18x and 32x** speedup. Similar results are observed for HITS and RWR as noted in Appendix F.

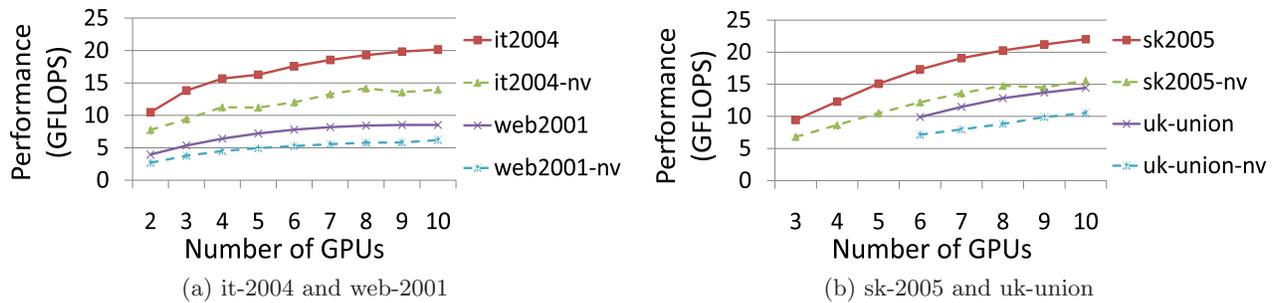

(a) it-2004 and web-2001

(b) sk-2005 and uk-union

Figure 4: Scalability of multi-GPU PageRank on web graphs and comparison with NVIDIA's HYB.

## 4.3 Multi-GPU PageRank on Web Graphs

Figure 4 shows the performance of computing PageRank with the multi-GPU framework introduced in Section 3.2 on the four web graph datasets in Table 3. The solid lines show the performance of using the Tile-Composite kernel; the dotted lines indicate the performance of using NVIDIA's HYB kernel. The lines for the sk-2005 and uk-union datasets start from 3 and 6 GPUs, because these datasets are very large and can only fit in the memory of at least 3 and 6 GPUs respectively. We observe that our multi-GPU framework and partition scheme can comfortably handle web graphs with billions of edges. For example on 10 GPUs, the distributed GPU implementation of PageRank can achieve about **23G-FLOPS** performance with **70%** parallel efficiency (on sk-2005) with the Tile-Composite kernel. On the two small datasets, it-2004 and web-2001, our Tile-Composite kernel achieves about **80%** parallel efficiency with 4 GPUs and **60%** parallel efficiency with 6 GPUs. All curves tend to flatten out after a point. This is because the workload size per GPU is low and the communication overheads begin to dominate and limit speedup. We should also emphasize that the performance of Tile-Composite kernel is about **1.55x** faster than HYB kernel on all datasets.

## 4.4 Parameter Auto-tuning

Next we will present experimental results to validate the auto-tuning method and the performance model introduced in Section 3.3. In these experiments we use the five matrices representing power-law graphs in Table 2. In our offline benchmarks, we set the upper bound of the workload size to 32768. This number is sufficiently large for most practical matrices that fit on the Tesla architecture, because there will be at least 960 warps (full occupancy) in each kernel, which correspond to about 31M non-zero entries in each tile.

In the first experiment we compare our heuristic approach for determining the number of tiles in our composite strategy. Figure 5(a) shows the number of tiles from the exhaustive search and the auto-tuning method. On the Webbase and Wikipedia matrices, our auto-tuned parameters are exactly the same as exhaustive searched results. We can see that our predicted number of tiles are very close to the optimal numbers on the other three datasets.

Figure 5(b) presents the optimal performance number by exhaustive search versus the performance number produced by using the auto-tuned number of tiles and partitioning strategies. The blue bars represent the optimal performance of exhaustive search; the yellow bars represent the results by running the tile-composite kernel with the auto-tuned parameters. On the Webbase and the Wikipedia matrices, the auto-tuning method achieves optimal performance. On the other datasets, the auto-tuned performance is within 3% of the optimal performance which is an excellent result.

Our auto-tuning method only needs the performance model to predict the *relative* performance trend under different parameter settings so as to automatically select the optimal parameter. However we also want to evaluate how accurately our performance model can predict the absolute performance. Figure 5(c) presents the measured and the predicted performance with the same parameters produced by auto-tuning. The blue bars are the results by running the kernel on the GPU, while the yellow bars are the simulation results using our performance model. We can see the predictions are accurate, and they are all within roughly 20% of the measured results. The error is largely attributable to the fact that we use the average performance of the different warps to estimate the overall performance, and the fact that the lookup table relies on synthetic benchmarks in which all workloads are of the same shape. To reiterate we should note that this error in prediction of actual performance does not significantly impact our auto-tuning method since their relative performance is what matters.

## 5. DISCUSSION

The power-law property is commonly observed in the datasets of graph mining problems. Our work is intended to find a better representation of matrices on GPU that is suitable for power-law graphs. Our extensive experiments have demonstrated the effectiveness of our optimizations on matrices representing power-law graphs. Next, we discuss the broader applicability of our optimization techniques.

**Tiling:** Our partially tiling optimization employs a greedy heuristic. The benefit of tiling comes from the fast accesses of the vector $x$ from the texture cache. The cost of each tile is the random writes to the resulting vector $y$. Both the benefit and the cost increase with the number of non-zeros in a tile. But the cost of random writes is bounded by the length of the vector $y$. After the cost of a tile reaches this bound, the more non-zero elements, the more benefit. Our tiling optimization starts with the densest columns, which is the most beneficial, and greedily finds the following tiles. If the non-zeros of a matrix concentrate in the first few tiles, our partially tiling optimization can efficiently finishes the majority of computation work without paying extra cost for the sparse columns. Power-law matrices are a subset of such matrices. Our experimental results are consistent with the above discussion. The only difference between COO and tile-coo kernel is tiling. On power-law matrices, tile-coo kernel performs consistently better than COO (Figure 2). On non-power-law matrices, tile-coo kernel is better than COO, but the benefit is very marginal (Figure 7 in Appendix D).

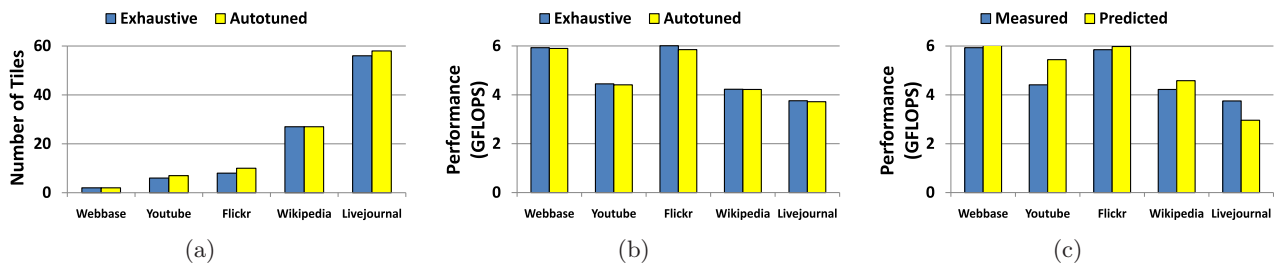

Figure 5: Performance auto-tuning and modeling. (a)Auto vs Exhaustive searched number of tiles. (b)Auto vs Exhaustive searched performance. (c)Predicted vs Measured performance of autotuned kernels.

**Composite Storage:** Our composite storage optimization consists of two parts: combination of CSR and ELL storage, and padding workload to warp size. The tile-composite kernel performs better than tile-coo kernel on both power-law and non-power-law matrices. This storage scheme can be applied to general matrices, although memory overhead of padded zeros should be considered as a constraint.

**Performance Modeling:** Our performance model improves the practical utility of the proposed optimizations. We only need to build the performance model once for the same hardware. The model does not rely on the power-law property of the matrix. Given an arbitrary matrix, we can get a relatively accurate prediction of the performance of our tile-composite kernel from the performance model before conducting large-scale experiments. More importantly, the CSR, CSR-vector and ELL kernels from NVIDIA can be modeled as special cases of our tile-composite kernel under the framework of our performance model. The CSR and CSR-vector kernel can be treated as the tile-composite kernel with a single tile and only CSR storage in the composite storage scheme; The ELL kernel can be seen as the tile-composite kernel with a single tile and with only ELL storage. With the generality of our performance model, the performance of different kernels can be predicted by plugging in the data to the model first. The best predicted kernel can be chosen to perform real computation of the data.

## 6. CONCLUSIONS

We proposed architecture conscious optimizations for the sparse matrix-vector multiply kernel on GPUs and studied the implications of this effort for graph mining. Our optimizations take into account both the architecture features of GPUs and the characteristics of graph mining applications. Our tiling approach utilizes the texture cache on GPUs more efficiently than previous work and provides much better memory locality. Our tiling with composite representation leverages the power-law characteristics of large graphs. We have obtained significant performance improvement over the state-of-the-art on such graphs. We also present empirical evaluations of applying our optimizations to PageRank, Random Walk with Restart and HITS algorithms. On these algorithms, our best kernel is **1.8 to 2.1** times faster than an industrial strength GPU competitor and from **18 to 32** times faster than a similar CPU implementation. The high performance of our optimizations relies on carefully tuning of parameters. We proposed a performance model to automatically tune our tile-composite kernel. We further extend our optimizations to handle web graphs on an MPI-based cluster.

**Acknowledgments:** This work is supported in part by grants from National Science Foundation CAREER-IIS-034-7662, RI-CNS-0403342, CCF-0702587 and IIS-0917070.

# APPENDIX

## A. GPU BACKGROUND

In this section, we discuss the hardware architecture and the programming model of CUDA GPUs. Figure 6 illustrates the organization of the computing hardware and the memory hierarchy in CUDA GPUs.

A CUDA device consists of a set of *streaming multiprocessors* (SMs), each one equipped with one instruction unit and 8 *streaming processors* (SPs). The parallel region of a CUDA program is partitioned into a *grid* of thread *blocks* that run *logically* in parallel. The programmer can decide the dimensions of the grid and the block. Thread blocks are distributed evenly on the multiprocessors. A *warp* is a group of 32 threads that run concurrently on a SM. The execution of the threads follows a *single instruction multiple threads* (SIMT) model. The instruction unit on a SM issues one instruction for all the threads in the same warp at each time [13]. The SPs executes this instruction for all the threads in the warp. Different warps within a block are time-shared on the hardware resources. A *kernel* is the code in the parallel region to be executed by each thread. Conditional instructions cause a divergence in the execution if threads in the same warp take different conditional paths. The threads are serialized in this situation.

There are various memory units on a CUDA device. The device memory, which is also called the global memory, is a large memory which is visible to all threads on the device [13]. The access latency of the global memory is high. Memory requests of a half warp (16 threads) are served together at a time. When accessing a 4- or 8-byte word, the global memory is organized into 128-byte segments [13]. The number of memory transactions executed for a half warp is the number of memory segments requested by this half warp.[1] The requests from the threads in the half warp are *coalesced* into one memory transaction if they are accessing addresses in the same segment. When the addresses accessed by the half warp are all in one segment, this request is *fully coalesced*. The global memory is divided into 8 equally-sized memory partitions of 256-byte width. Concurrent memory requests to the global memory by all the active warps should be distributed uniformly amongst partitions. The term *partition camping* [14] is used to describe the situation when global memory accesses are congested and queued up at some partitions while the other partitions are idle. While *coalescing* concerns global memory accesses within a half warp, *partition camping* concerns global memory accesses amongst all active half warps. The constant and texture memories are read-only regions in the global memory space with on-chip caches. The programmer can bind a region of the global memory to either the constant or the texture memory before the kernel starts. Each multiprocessor is equipped with an on-chip scratchpad memory [13], which is called the shared memory. The shared memory has very low access latency. It is only visible to the threads within one block and has the same lifetime as the block [13]. The shared memory is organized into banks. If multiple addresses in the same bank are accessed at the same time, it leads to bank conflicts and the accesses are serialized. There are also a set of registers shared by the threads in the block.

---

[1] Devices with Compute Capability lower than 1.2 have stricter requirements.

| Matrix | Rows | Columns | NNZ | Power-law? |
|---|---|---|---|---|
| Dense | 2K | 2K | 4M | No |
| Circuit | 171K | 171K | 0.96M | No |
| FEM/Harbor | 47K | 47K | 2.4M | No |
| LP | 4.3K | 1M | 11M | No |
| Protein | 36K | 36K | 4M | No |
| Webbase | 1M | 1M | 3M | Yes |
| Flickr | 1.7M | 1.7M | 22.6M | Yes |
| LiveJournal | 5.2M | 5.2M | 77M | Yes |
| Wikipedia | 1.9M | 1.9M | 40M | Yes |
| Youtube | 1.1M | 1.1M | 4.9M | Yes |

Table 2: Matrix and Graph Datasets

## B. NVIDIA'S SPMV LIBRARY

The Sparse Matrix-Vector Multiplication (SpMV) kernel computes a vector $y$ as the product of a $n$ by $m$ sparse matrix $A$ and a dense vector $x$. The *compressed sparse row* (CSR) format stores non-zeros in the same row contiguously in memory, and all rows are stored in one data array, with another array holding the column indices of the non-zeros. A third array of row pointers marks the boundary of each row. The corresponding CSR kernel [3] assigns the computation of each row to a thread. With power-law graphs, it is hard to balance the workload among threads within one thread block. So all the threads in one block will wait for the thread which is assigned to the longest row. To improve this method, CSR-vector format [3] uses a warp of 32 threads to work on each row. This strategy only helps the rows with more than 32 non-zeros, but most of the nodes in power-law graphs have degree lower than 32. The computation resources of the warps assigned to such rows will therefore be wasted. Baskaran and Bordawekar [2] further optimized the CSR-vector format by using a half warp for each row to improve global memory accesses, and also a padding technique is used to ensure the memory requests are fully coalesced. But there are still considerable amount of threads wasted on the rows with less than half a warp of non-zero elements.

Besides the CSR format, the *coordinate* (COO) and *ellpack* (ELL) formats are also used in Bell and Garland's SpMV kernel. In COO format, all the non-zeros in matrix $A$ are combined into a long vector grouped by row index, and the kernel first computes the multiplication of each non-zeros with the corresponding elements of vector $x$ in the first pass; then the segmented reduction of the rows is done on this long vector by thread warps. In the reduction phase, because the length of each row is not necessarily a multiple of warp size, synchronization points are heavily used and warp thread divergence is frequent. However, the COO kernel is the most insensitive to variable row length in the matrix according to the previous study [3]. The ELL format requires the number of non-zeros on each row is bounded by some small number $k$, so that the matrix $A$ can be represented by a dense $n$ by $k$ matrix $M$, in which only non-zeros in $A$ are stored, and the corresponding column indices of these non-zeros are also stored in a separate matrix. In the ELL kernel, $M$ is stored in column major, and the thread assigned to each row can access global memory very efficiently. Zeros are added to rows with fewer than $k$ non-zeros, so $k$ cannot be large, otherwise it will introduce large overhead to access these zeros. The ELL format cannot be directly applied to graph mining algorithms, where the node degree in the graph cannot be bounded by a small number $k$. However, ELL and COO formats can be mixed together to represent a matrix, where the first $k$ non-zeros of each row are stored in ELL format and the others are stored in COO format. This is the *hybrid* (HYB) kernel of NVIDIA's SpMV Library [3]. There are two other formats in NVIDIA's SpMV Library [3]. The *diagonal* (DIA) format is only applicable to matrices in which all non-zeros fall into a band around the diagonal. The *packet* (PKT) format first uses Metis [9] to cluster non-zeros into dense sub-blocks, then a sub-block is loaded into shared memory and processed by a thread block as a dense sub-matrix. But the code of these two kernels cannot run on matrices of power-law graphs in our experiments.

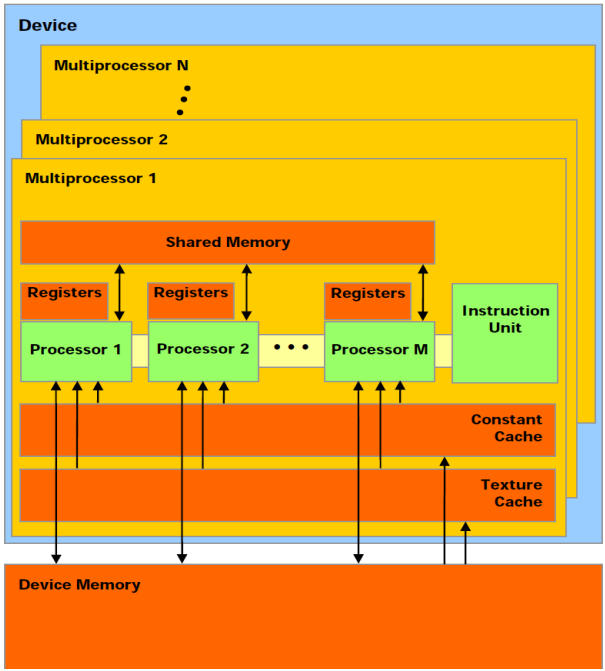

Figure 6: Hardware Organization and Memory Hierarchy of a CUDA Device [14]

## C. DATASET AND HARDWARE DETAIL

**Datasets**: In our single GPU experiments, we use four web-based graph datasets, including user link relationship graphs from Flickr, LiveJournal and Youtube and a web-page link relationship graph from Wikipedia [11]. All graphs exhibit power-law characteristics. In addition to the graph datasets we also include results on six popular unstructured matrix datasets used in previous studies [3]. Among these, one is a 2000 by 2000 dense matrix, which while not sparse, is a useful benchmark to show the maximum bandwidth that each kernel can achieve. Details of these graphs (represented in an adjacency matrix) and matrices are shown in Table 2. In the four graph datasets, the number of non-zeros (NNZ) is the number of directed links and the number of rows (or columns) is the number of nodes in the graphs. In our multi-GPU experiments, the web graph datasets used are provided in Table 3. These web graphs were crawled using *UbiCrawler* [4] developed by the *Laboratory for Web Algorithmics* at the *Univerita Degli Studi Di Milano*. None of these datasets can fit in the memory of one GPU.

**Hardware configuration**: Our experiments are run in an MPI-based cluster environment. On the CPU side, each node has an Opteron X2 2218 CPU with 8 GB of main memory. On the GPU side, each node is equipped with two

| Graph | Nodes | Edges | Power-law? |
|---|---|---|---|
| it-2004 | 41,291,594 | 1,150,725,436 | Yes |
| sk-2005 | 50,636,154 | 1,949,412,601 | Yes |
| uk-union | 133,633,040 | 5,507,679,822 | Yes |
| web-2001 | 118,142,155 | 1,019,903,190 | Yes |

Table 3: Web Graph Datasets

NVIDIA Tesla C1060 GPUs. Each GPU has 30 multiprocessors with 240 processing cores and 4 GB of global memory. The single GPU experiments are run on a single node and a single GPU. The multi-GPU experiments are run on multiple nodes, and each node uses a single GPU. The CPU code is complied with the gcc compiler version 4.1.2. The GPU code is compiled with CUDA version 3.0.

## D. PERFORMANCE ON UNSTRUCTURED MATRIX AND COMPARISON WITH CPU

**Performance on Unstructured Matrix Data**: The speed and bandwidth performance of different kernels on non-power-law matrices are shown in Figure 7. We immediately observe that our methods while comparing favorably on some of the kernels do not always perform as strongly as the best. In fact on these datasets, interestingly, there is no single kernel that outperforms all others.

Our tiling with composite storage kernel performs the best on the 2000 by 2000 dense matrix with 17.57 GFLOPS speed and 105.5 GB/s bandwidth. This bandwidth utilization is higher than the peak bandwidth of 102 GB/s in the official hardware specification from NVIDIA website. This somewhat surprising result is due to the effect of texture binding of vector $x$ allowing for elements in $x$ to be directly fetched from the cache. Our tiling with composite storage kernel runs 30% faster than CSR-vector kernel on the dense matrix. This is because we pad the storage of the matrix in global memory to ensure that all global memory accesses are fully coalesced. The CSR-vector format concatenates all rows together. If one row is not padded to an integer multiple of the warp size, all global memory accesses after this row will not be fully coalesced resulting in a loss in performance.

Baskaran and Bordawekar's CSR kernel performs best on FEM/Harbor and Protein dataset. Their kernel achieves 12.76 GFLOPS speed and 78.6 GB/s bandwidth on FEM/Harbor, and 15.74 GFLOPS speed and 95.5 GB/s bandwidth on Protein. The bandwidth utilizations are close to the maximum on these two datasets. HYB kernel performs best on the other two non-power-law matrices. It achieves 5.98 GFLOPS speed and 45.4 GB/s bandwidth on Circuit, and 8.45 GFLOPS speed and 61.6 GB/s bandwidth on LP. On all four non-power-law matrices, our tiling with composite storage kernel is amongst the top four in speed and bandwidth. Our tiling with composite storage kernel is only 10.5% slower than HYB kernel on Circuit matrix.

The non-zero elements present a relatively balanced distribution in these non-power-law matrices, not as biased as power-law matrices. Our methods first reorder the columns and partition the matrices into a denser and a sparser submatrix. This partition will not produce a dense enough matrix in which most of the non-zeros in the original matrix are concentrated. This phenomenon leads to the low performance of our methods on non-power-law matrices, because the tiling of the denser matrix still requires overhead, but does not gain benefit in performance.

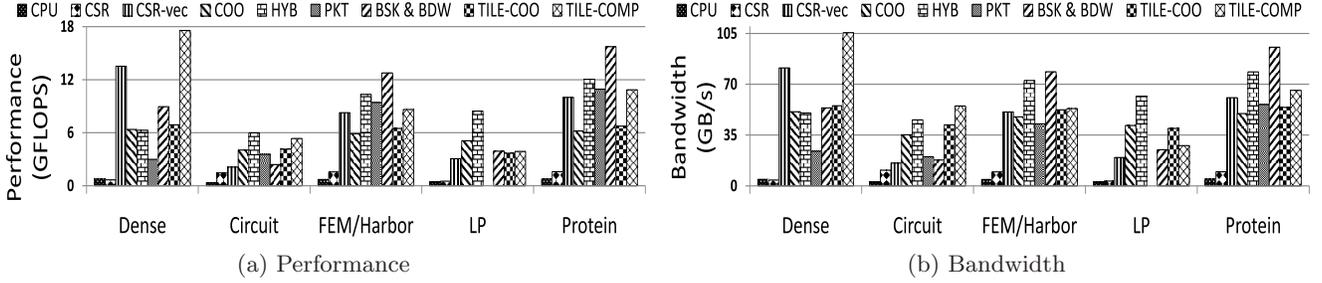

(a) Performance  (b) Bandwidth

Figure 7: SpMV kernels comparison on unstructured matrices from NVIDIA's SpMV Library [3].

**Comparison with CPU SpMV**: Previous works [3, 6] have already illustrated the benefits of GPU vs CPU. Our main point is to demonstrate the benefits of our approach over other GPU work on power-law graphs. The CPU results are included for the completeness of our evaluation. We implemented CSR SpMV kernel on the CPU. CSR format is the most efficient on CPU among different sparse matrix formats. We ran experiments with the CPU kernel on all datasets in Table 2. The GPU kernels significantly outperform the CPU kernel in almost all settings. GPU CSR kernel is the slowest kernel on GPU. It is slower than CPU kernel on the Dense matrix data because the clock rate of one GPU processor is lower than CPU. The GPU kernels perform dominantly faster than CPU kernel in all the other formats with speedups ranging from 2.05x to 37.31x.

## E. AUTOTUNING PSEUDOCODE

**Algorithm 1** Tile-Composite Kernel Auto-tuning

1: **Input:** $n$ by $n$ matrix $M$ sorted by column lengths
2: **Output:** number of tiles and partition size of each tile
3: $TILE\_WIDTH \leftarrow 64K; NTile \leftarrow 0$
4: **while** $NTile < n/TILE\_WIDTH$ **do**
5: $\quad StartCol \leftarrow NTILE \times TILE\_WIDTH$
6: $\quad$ **if** $M.ColLength[StartCol] \leq 1$ **then**
7: $\quad\quad$ break;
8: $\quad$ **else**
9: $\quad\quad$ **for** $i = StartCol$ to $StartCol + TILE\_WIDTH$ **do**
10: $\quad\quad\quad InsertCol(M.Tile[NTile], M.Col[i])$
11: $\quad\quad$ **end for**
12: $\quad\quad WL = Partition(M.Tile[NTile])$
13: $\quad\quad NTile \leftarrow NTile + 1$
14: $\quad$ **end if**
15: **end while**
16: **Return:** $NTile, WL$

**Algorithm 2** $Partition(T)$: Partition of one tile of matrix

1: **Input:** one tile $T$ from the matrix sorted by row lengths
2: **Output:** optimal workload size to partition $T$
3: $WL_{low} \leftarrow T.RowLength[0]$ {Workload lower bound}
4: $WL_{up} \leftarrow \frac{T.NNZ}{MAX\_ACT\_WARP}$ {Workload upper bound}
5: $OptWL \leftarrow 0; OptTime \leftarrow +\infty; WL \leftarrow WL_{low}$
6: **while** $WL \leq WL_{up}$ **do**
7: $\quad Time = PM(T, WL)$ {Performance Modeling of T with WL}
8: $\quad$ **if** $Time < OptTime$ **then**
9: $\quad\quad OptTime \leftarrow Time; OptWL \leftarrow WL$
10: $\quad$ **end if**
11: $\quad WL \leftarrow WL + T.RowLength[0]$
12: **end while**
13: **Return:** $OptWL, OptTime$

## F. GRAPH MINING ALGORITHMS

A large class of graph mining algorithms leverage the SpMV kernel iteratively to perform computation until the algorithms converge, e.g. PageRank [15], HITS [10] and

**Algorithm 3** $PM(T, WL)$: Performance Modeling of tile T given workload size

1: **Input:** tile T and workload size WL
2: **Output:** Total run time
3: $NWarp \leftarrow \lceil \frac{T.NNZ}{WL} \rceil$
4: $MAX\_ACT\_WARP \leftarrow MAX\_ACT\_WARP/SM * NUM\_SM$
5: $I = \lceil \frac{NWarp}{MAX\_ACT\_WARP} \rceil$ {Number of iterations}
6: $i \leftarrow 0; j \leftarrow 0$ {Row Index; Warp Index}
7: **while** $i < T.NumRow$ **do**
8: $\quad$ {Partition T with $workload = WL$}
9: $\quad w_j \leftarrow T.RowLength[i]; h_j \leftarrow \frac{WL}{w_j}$
10: $\quad Padding(w_j, h_j, WarpSize)$ {Padding w or h}
11: $\quad IterId = \frac{j}{MAX\_ACT\_WARP}$
12: $\quad P[IterId] \leftarrow P[IterId] + Performance(w_j, h_j)$
13: $\quad Size[IterId] \leftarrow Size[IterId] + w_j \times h_j$
14: $\quad Cnt[IterId] \leftarrow Cnt[IterId] + 1$
15: $\quad j \leftarrow j + 1, i \leftarrow i + h_j$
16: **end while**
17: **for** $i = 0$ to $I - 1$ **do**
18: $\quad P[i] \leftarrow \frac{P[i]}{Cnt[i]}; t[i] \leftarrow \frac{Size[i]}{P[i]}$
19: $\quad TotalTime \leftarrow TotalTime + t[i]$
20: **end for**
21: **Return:** $TotalTime$

Random Walk with Restart [18]. These algorithms first transform the adjacency matrix of a graph and then operate on the transformed matrix. The graph dataset used by these algorithms usually have strong power-law properties, hence the number of non-zeros on each row or column of the corresponding matrix will follow a power-law distribution. The skew of the distribution leads to poor load balancing and low memory access efficiency on GPU.

In this section, we describe three data mining algorithms which can be written in the form of matrix-vector multiplication. These algorithms essentially compute the power method for different matrices related to the link structure of the graphs. Within one iteration of the power method, the running time is dominated by the time required to compute the matrix-vector product. These algorithms usually operate on large power-law graphs. Hence they can be sped up using our sparse matrix representation and computed by our SpMV kernels.

**PageRank**: The PageRank algorithm models the link structure of web graphs by the random walk behavior of a *random surfer* [15]. The web graph can be represented by a directed graph $G = (V, E)$, where $V$ is a set of $n$ vertices and $E$ is the set of directed edges. The adjacency matrix $A$ is defined as $A(u,v) = 1$ if edge $(u,v) \in E$; otherwise, $A(u,v) = 0$. Matrix $W$ denotes the row normalized matrix of $A$. The PageRank vector $p$ is computed iteratively using the following equation until it converges:

$$p^{(k+1)} = cW^T p^{(k)} + (1-c)p^{(0)} \qquad (6)$$

241

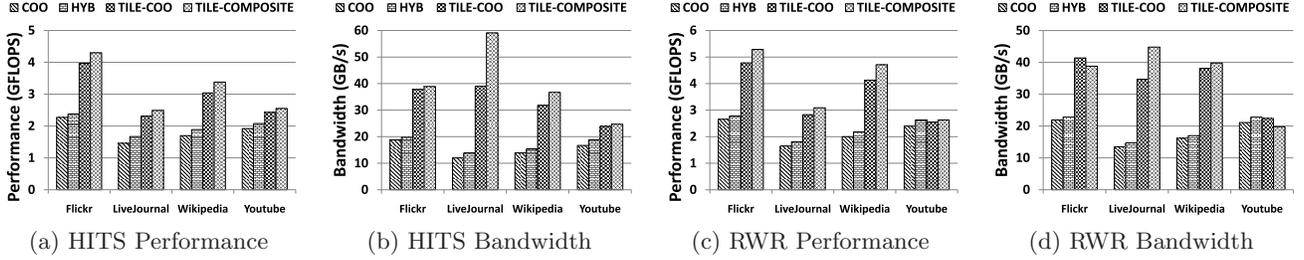

Figure 8: Performance and bandwidth of HITS and RWR on graph datasets.

| Graph | CPU | COO | HYB | TILE-COO | TILE-Comp |
|---|---|---|---|---|---|
| Flickr | 4.97 | 0.40 | 0.38 | 0.23 | 0.21 |
| LiveJournal | 44.88 | 3.82 | 3.33 | 2.41 | 2.24 |
| Wikipedia | 39.36 | 2.73 | 2.45 | 1.52 | 1.37 |
| Youtube | 4.35 | 0.33 | 0.30 | 0.26 | 0.25 |

Table 4: Running time of HITS (in seconds)

| Graph | CPU | COO | HYB | TILE-COO | TILE-Comp |
|---|---|---|---|---|---|
| Flickr | 8.25 | 0.59 | 0.56 | 0.33 | 0.29 |
| LiveJournal | 36.99 | 2.85 | 2.60 | 1.73 | 1.52 |
| Wikipedia | 23.23 | 1.46 | 1.35 | 0.71 | 0.62 |
| Youtube | 2.32 | 0.14 | 0.13 | 0.14 | 0.13 |

Table 5: Average running time of Random Walk with Restart (in seconds) on 25 random query nodes

where $c$ is a damping factor (set to 0.85), $p^{(0)}$ is initialized as a $n$ by 1 vector with all elements set to $1/n$.

**HITS**: HITS is a link analysis algorithm of web pages [10]. It gives each web page two attributes: authority and hub. It rates web pages by assigning authority score and hub score to each web page. Let matrix $A$ be the adjacency matrix of a directed graph $G = (V, E)$ or $G$ may be a query specific subgraph of the whole web graph. Then the authority score vector $\vec{a}$ and hub score vector $\vec{h}$ are recursively defined as

$$\vec{a}^{(k+1)} = A^T \vec{h}^{(k)} \qquad \vec{h}^{(k+1)} = A \vec{a}^{(k)} \qquad (7)$$

This recursive definition with two matrix-vector products can be rewritten as one matrix and vector multiplication by

$$\left[ \begin{array}{c} \vec{a} \\ \vec{h} \end{array} \right]^{(k+1)} = \left[ \begin{array}{cc} 0 & A^T \\ A & 0 \end{array} \right] \left[ \begin{array}{c} \vec{a} \\ \vec{h} \end{array} \right]^{(k)} \qquad (8)$$

The power method can be used to solve this eigen vector problem. Elements in $\vec{a}^{(0)}$ and $\vec{h}^{(0)}$ vectors are all initialized to $1/|V|$. In each iteration, a $2|V|$ by $2|V|$ matrix in Equation 8 is multiplied by a vector combined with $\vec{a}$ and $\vec{h}$. Then the first and second half of the resulting vector are normalized to sum to 1 separately. Each normalization requires a reduction operation on the vector and a division of the vector by a constant. A convergence check is also needed at the end of each iteration. Each iteration of our HITS implementation involves one SpMV kernel, three parallel reduction kernels (two for normalization and one for convergence check) and two vector division by constant kernels. The vector division by constant kernel can be implemented very efficiently in the same way as vector addition. On our implementation of the HITS algorithm we compare the performance of our four GPU SpMV kernels on the four graph datasets. The speed and bandwidth performance are shown in Figure 8(a) and Figure 8(b). Our TILE-COO and TILE-Composite kernels perform better than COO and HYB kernels in all four datasets. On Flickr, LiveJournal and Wikipedia, the speedups are similar to those observed in PageRank algorithm. On Youtube, our optimizations are actually a bit faster when compared to the NVIDIA kernels in spite of the relatively small size of the dataset. Combining the two matrices into one in the HITS algorithm results in a larger and sparser matrix making it more amenable to our optimizations. The total running time compared with CPU implementation is listed in Table 4. We observe a 17x to 29x speedup of the GPU implementations over the corresponding CPU implementation.

**Random Walk with Restart**: Random Walk with Restart (RWR) is an algorithm that tries to measure the relevance between two nodes in an undirected graph [18]. Given a query node $i$ in the graph, the relevance score from all other nodes to node $i$ forms a vector $\vec{r_i}$. In RWR, vector $\vec{r_i}$ is computed by the following equation:

$$\vec{r_i}^{(k+1)} = cW\vec{r_i}^{(k)} + (1-c)\vec{e_i} \qquad (9)$$

where $c$ is a restart probability parameter (set to 0.9 in our experiment), $W$ is the column normalized adjacency matrix and $\vec{e_i}$ is a vector whose $i^{th}$ element is 1 and all the other elements are 0. Vector $\vec{r_i}$ can be computed using the power method. In each iteration, there is a matrix-vector multiplication followed by a vector addition and a convergence checking operation. In our implementation, we use the GPU SpMV kernels for matrix-vector multiplication, and GPU parallel reduction for checking convergence in the same way as PageRank. An efficient vector addition kernel is also implemented by assigning one GPU thread to compute one element in the resulting vector. Note that RWR is an interactive application, we randomly select 25 query nodes and the performance is reported by averaging (arithmetic mean) the result of each query. The number of computations per iteration is the same whichever node is selected as query, so the experiment results of the randomly selected 25 query nodes can reflect the speed of different kernels. Since RWR operates on undirected graphs, we treat each link in our directed graph datasets as an undirected link. The speed and bandwidth performance of RWR implementations on four graph datasets based on four GPU SpMV kernels are shown in Figure 8(c) and Figure 8(d). The total running time is listed in Table 5. We observe similar performance results as in the case of PageRank. Our optimized TILE-COO and Tile-Composite kernels are 1.5x to 2.0x as fast as COO and HYB kernels on Flickr, LiveJournal and Wikipedia graphs. The four kernels perform about the same on Youtube graph. All GPU implementations are 13x to 37x faster than CPU implementation. The best speedup is achieved by our TILE-Composite kernel on Wikipedia graph.